\documentclass[12pt]{amsart}

\usepackage{enumitem}
\usepackage{psfrag}
\usepackage{graphicx}
\usepackage{pinlabel}
\usepackage{graphicx}
\usepackage{amsmath}
\usepackage{amssymb}
\usepackage{amscd}
\usepackage{autobreak}
\usepackage{picinpar}
\usepackage{times}
\usepackage{pb-diagram}
\usepackage{graphicx}
\usepackage{wrapfig}
\usepackage{xspace}
\usepackage{hyperref}
\usepackage{url}
\usepackage{caption}
\usepackage{subcaption}
\usepackage{fancyhdr}
\usepackage{overpic}
\usepackage{color}
\usepackage[square,comma,sort&compress,numbers]{natbib}
\usepackage{latexsym}
\usepackage{mathrsfs}
\usepackage{amsfonts}
\usepackage{hyperref}
\usepackage{amsthm}
\usepackage{appendix}
\usepackage{pdfsync}

\theoremstyle{definition}
\newtheorem{theorem}{Theorem}[section]
\newtheorem{lemma}[theorem]{Lemma}
\newtheorem{proposition}[theorem]{Proposition}
\newtheorem{corollary}[theorem]{Corollary}
\newtheorem{definition}[theorem]{Definition}

\newtheorem{remark}[theorem]{Remark}
\newtheorem*{theorem*}{Theorem}
\def\qed{\hfill{Q.E.D.}\smallskip}

\begin{document}

\title{\bf Discrete conformal structures on surfaces with boundary (III)---Deformation }
\author{Xu Xu, Chao Zheng}

\date{\today}

\address{School of Mathematics and Statistics, Wuhan University, Wuhan, 430072, P.R.China}
 \email{xuxu2@whu.edu.cn}

\address{International Center for Mathematical Research (BICMR), Beijing (Peking) University, Beijing, 100871, P.R. China}
\email{czheng@whu.edu.cn}



\begin{abstract}
The present work constitutes the third installment in a series of investigations devoted to discrete conformal structures on surfaces with boundary.
In our preceding works \cite{X-Z DCS1, X-Z DCS2}, we established, respectively, a classification of these discrete conformal structures and results on their rigidity and existence.
Building on this foundation, the present work focuses on the deformation theory of discrete conformal structures on surfaces with boundary.
Specifically, we introduce the combinatorial Ricci flow and the combinatorial Calabi flow,
and establish the longtime existence and global convergence of solutions to these combinatorial curvature flows.
These results yield effective algorithms for finding discrete hyperbolic metrics on surfaces with totally geodesic boundaries of prescribed lengths.
\end{abstract}

\maketitle

\section{Introduction}

\subsection{Background}

In discrete conformal geometry, a central problem is finding discrete conformal structures with prescribed combinatorial curvatures.
The combinatorial curvature flow, inspired by Chow-Luo's foundational work on combinatorial Ricci flow for Thurston's circle packings on polyhedral surfaces \cite{Chow-Luo}, has emerged as an effective approach to address this problem.
Following this line of inquiry, Ge \cite{Ge1,Ge2} extended the framework by introducing the combinatorial Calabi flow for Thurston's circle packings on closed surfaces.
While extensive research has been devoted to exploring combinatorial curvature flows for discrete conformal structures on closed surfaces, investigations into their counterparts for surfaces with boundary remain relatively limited.

Motivated by Thurston's theory of circle packings on closed surfaces \cite{Thurston},
Guo-Luo \cite{GL2} pioneered the study of generalized circle packings on surfaces with boundary and subsequently proposed a generalized combinatorial curvature flow.
In our previous work \cite{XZ}, we established the longtime existence and global convergence of solutions to this generalized combinatorial curvature flow for a specific class of generalized circle packings, specifically the $(-1,-1,-1)$-type.
Furthermore, we introduced the corresponding combinatorial Calabi flow and proved the longtime existence and global convergence of its solutions.
Building upon Luo's work on vertex scaling of piecewise linear metrics on closed surfaces \cite{Luo1},
Guo \cite{Guo} introduced a class of hyperbolic discrete conformal structures, referred to as vertex scalings, on surfaces with boundary.
Additionally, a combinatorial Yamabe flow was constructed in this framework.
Motivated by Guo's combinatorial Yamabe flow, Li-Xu-Zhou \cite{Li-Xu-Zhou} recently introduced a modified combinatorial Yamabe flow for Guo's vertex scaling on surfaces with boundary,
which generalizes and completes Guo's results.
Additionally, inspired by \cite{Ge1,Ge2,GX3,Wu-Xu}, Luo-Xu \cite{Luo-Xu} introduced the combinatorial Calabi flow and the fractional combinatorial Calabi flow for Guo's vertex scaling on surfaces with boundary.
Recently, Xu \cite{Xu22} introduced a new class of hyperbolic discrete conformal structures on surfaces with boundary, together with the corresponding combinatorial Ricci flow, combinatorial Calabi flow, and fractional combinatorial Calabi flow in this framework.

In the present work, building upon the aforementioned studies, we introduce combinatorial curvature flows, specifically the combinatorial Ricci flow and the combinatorial Calabi flow, for discrete conformal structures on surfaces with boundary,
which unifies and generalizes the existing combinatorial curvature flows for discrete conformal structures on surfaces with boundary.
Moreover, we establish results on the longtime existence and global convergence of solutions to these combinatorial curvature flows on surfaces with boundary.

\subsection{Basic setting and previous results}

Let $\widetilde{\mathcal{T}}=(V,\widetilde{E},\widetilde{F})$ be a triangulation of a closed surface $\widetilde{\Sigma}$,
where $V$, $\widetilde{E}$, and $\widetilde{F}$ denote the sets of vertices, unoriented edges, and faces, respectively.
Here $V$ is a finite subset of $\widetilde{\Sigma}$ with $|V|=N$.
Let $N(V)$ be a small open regular neighborhood of the union of all vertices, chosen such that the neighborhoods around distinct vertices are disjoint.
Then $\Sigma=\widetilde{\Sigma}\setminus N(V)$ is a compact surface with $N$ boundary components.
The intersection $\mathcal{T}=\widetilde{\mathcal{T}}\cap \Sigma$ is referred to as an ideal triangulation of $\Sigma$.
Let $E=\widetilde{E}\cap \Sigma$ and $F=\widetilde{F}\cap \Sigma$ denote the sets of unoriented ideal edges and ideal faces of $\Sigma$ with respect to $\mathcal{T}$, respectively.
A boundary arc is defined as a connected component of the intersection between an ideal face and the boundary $\partial\Sigma$.
For simplicity, we index the boundary components of $\Sigma$ by the set $B = \{1, 2, \dots, N\}$, where each element $i \in B$ corresponds to a unique boundary component.
Let $E_+$ denote the set of oriented ideal edges.
For any two adjacent boundary components $i, j \in B$, the unoriented ideal edge connecting them is denoted by $\{ij\} \in E$, with the corresponding oriented ideal edge denoted by $(i, j)$.
An ideal face $\{ijk\} \in F$ is adjacent to the boundary components $i, j, k \in B$.
Specifically, such an ideal face is a hexagon, which corresponds to the triangle $v_iv_jv_k\in \widetilde{F}$ in $\widetilde{\mathcal{T}}$.
The sets of real-valued functions on $B$, $E$, and $E_+$ are denoted by $\mathbb{R}^N$, $\mathbb{R}^E$, and $\mathbb{R}^{E_+}$, respectively.

A edge length function associated to $\mathcal{T}$ is a vector $l \in \mathbb{R}^E_{>0}$ that assigns to each ideal edge $\{ij\} \in E$ a positive number $l_{ij} = l_{ji}$.
For any ideal face $\{ijk\} \in F$, there exists a unique right-angled hyperbolic hexagon (up to isometry) whose three non-adjacent edges have lengths $l_{ij}, l_{ik}, l_{jk}$ (see \cite[Theorem 3.5.14]{Ratcliffe}).
By gluing all such right-angled hyperbolic hexagons isometrically along the ideal edges in pairs,
one can construct a hyperbolic surface with totally geodesic boundary from the ideal triangulation $\mathcal{T}$.
Conversely, any ideally triangulated hyperbolic surface $(\Sigma, \mathcal{T})$ with totally geodesic boundary induces a function $l \in \mathbb{R}^E_{>0}$, where $l_{ij}$ is the length of the shortest geodesic connecting the boundary components $i, j \in B$. Such an edge length function $l \in \mathbb{R}^E_{>0}$ is referred to as a \textit{discrete hyperbolic metric} on $(\Sigma, \mathcal{T})$.
The length $K_i$ of the boundary component $i \in B$ is called the \textit{generalized combinatorial curvature} of the discrete hyperbolic metric $l$ at $i$. Specifically, the generalized combinatorial curvature $K_i$ is defined by
\begin{equation}\label{Eq: K}
K_i = \sum_{\{ijk\} \in F} \theta^{jk}_i,
\end{equation}
where the summation is taken over all right-angled hyperbolic hexagons adjacent to $i$, and the generalized angle $\theta^{jk}_i$ is the length of the boundary arc of the right-angled hyperbolic hexagon $\{ijk\}$ at $i$.

The main result in \cite{X-Z DCS1} is the following theorem, which provides explicit forms and a classification of the discrete conformal structures on surfaces with boundary.

\begin{theorem}[\cite{X-Z DCS1}, Theorem 1.4]\label{Thm: DCS}
Let $(\Sigma,\mathcal{T})$ be an ideally triangulated surface with boundary, and let $d=d(f)$ be a discrete conformal structure on $(\Sigma,\mathcal{T})$.
There exist constant vectors $\alpha\in \mathbb{R}^N$ and $\eta\in \mathbb{R}^{E}$ satisfying $\eta_{ij}=\eta_{ji}$, such that for any right-angled hyperbolic hexagon $\{ijk\}\in F$,
the value of $\frac{\sinh d_{ij}}{\sinh d_{ji}}$ can only take one of the following four forms:
\begin{gather*}
\sqrt{\frac{1+\alpha_ie^{2f_i}}{1+\alpha_je^{2f_j}}},\
-\sqrt{\frac{1+\alpha_ie^{2f_i}}{1+\alpha_je^{2f_j}}},\
e^{\frac{1}{2}C_{ij}(f_i-f_j)},\
-e^{\frac{1}{2}C_{ij}(f_i-f_j)},
\end{gather*}
where $C\in \mathbb{R}^{E_+}$ is a constant vector satisfying
$C_{ij}+C_{jk}+C_{ki}=0$ for any $\{ijk\}\in F$ and $C_{rs}+C_{sr}=0$ for any subset $\{r,s\}\subseteq\{i,j,k\}$.
Furthermore,
\begin{description}
\item[(Ai)]
if $\frac{\sinh d_{ij}}{\sinh d_{ji}}=\sqrt{\frac{1+\alpha_ie^{2f_i}}{1+\alpha_je^{2f_j}}}>0$, then the discrete conformal structure $d=d(f)$ takes the following form
\begin{equation}\label{Eq: d3}
\coth d_{ij}
=-\frac{\alpha_ie^{2f_i}}{\sinh l_{ij}}\sqrt{\frac{1+\alpha_je^{2f_j}}{1+\alpha_ie^{2f_i}}}
+\frac{\eta_{ij}e^{f_i+f_j}}{\sinh l_{ij}}
\end{equation}
with
\begin{equation}\label{Eq: DCS3}
\cosh l_{ij}
=-\sqrt{(1+\alpha_ie^{2f_i})(1+\alpha_je^{2f_j})}
+\eta_{ij}e^{f_i+f_j},\ \text{for}\ 1+\alpha_ie^{2f_i}>0, 1+\alpha_je^{2f_j}>0,
\end{equation}
\begin{equation}\label{Eq: new 1}
\cosh l_{ij}
=\sqrt{(1+\alpha_ie^{2f_i})(1+\alpha_je^{2f_j})}
+\eta_{ij}e^{f_i+f_j},\ \text{for}\ 1+\alpha_ie^{2f_i}<0, 1+\alpha_je^{2f_j}<0;
\end{equation}

\item[(Aii)]
if $\frac{\sinh d_{ij}}{\sinh d_{ji}}=e^{\frac{1}{2}C_{ij}(f_i-f_j)}>0$, then the discrete conformal structure $d=d(f)$ takes the following form
\begin{equation}\label{Eq: d1}
\coth d_{ij}=\frac{\sinh(f_j-f_i-C_{ij})}{\sinh l_{ij}}+\frac{\eta_{ij}e^{f_i+f_j}}{\sinh l_{ij}}
\end{equation}
with
\begin{equation}\label{Eq: DCS1}
\cosh l_{ij}
=-\cosh(f_j-f_i-C_{ij})+\eta_{ij}e^{f_i+f_j};
\end{equation}

\item[(Bi)]
if $\frac{\sinh d_{ij}}{\sinh d_{ji}}=-\sqrt{\frac{1+\alpha_ie^{2f_i}}{1+\alpha_je^{2f_j}}}<0$, then the discrete conformal structure $d=d(f)$ takes the following form
\begin{equation*}
\coth d_{ij}
=\frac{\alpha_ie^{2f_i}}{\sinh l_{ij}}
\sqrt{\frac{1+\alpha_je^{2f_j}}{1+\alpha_ie^{2f_i}}}
+\frac{\eta_{ij}e^{f_i+f_j}}{\sinh l_{ij}}
\end{equation*}
with
\begin{equation}\label{Eq: DCS4}
\cosh l_{ij}
=\sqrt{(1+\alpha_ie^{2f_i})(1+\alpha_je^{2f_j})}
+\eta_{ij}e^{f_i+f_j},\ \text{for}\ 1+\alpha_ie^{2f_i}>0, 1+\alpha_je^{2f_j}>0,
\end{equation}
\begin{equation}\label{Eq: new 2}
\cosh l_{ij}
=-\sqrt{(1+\alpha_ie^{2f_i})(1+\alpha_je^{2f_j})}
+\eta_{ij}e^{f_i+f_j},\ \text{for}\ 1+\alpha_ie^{2f_i}<0, 1+\alpha_je^{2f_j}<0;
\end{equation}

\item[(Bii)]
if $\frac{\sinh d_{ij}}{\sinh d_{ji}}=-e^{\frac{1}{2}C_{ij}(f_i-f_j)}<0$, then the discrete conformal structure $d=d(f)$ takes the following form
\begin{equation*}
\coth d_{ij}=-\frac{\sinh(f_j-f_i-C_{ij})}{\sinh l_{ij}}+\frac{\eta_{ij}e^{f_i+f_j}}{\sinh l_{ij}}
\end{equation*}
with
\begin{equation}\label{Eq: DCS2}
\cosh l_{ij}
=\cosh(f_j-f_i-C_{ij})+\eta_{ij}e^{f_i+f_j}.
\end{equation}
\end{description}
Moreover, there exist precisely six distinct types of combinations of such discrete conformal structures on the surface with boundary, specifically: (\ref{Eq: DCS3}), (\ref{Eq: new 1}), (\ref{Eq: DCS1}), the combination of (\ref{Eq: DCS3}) and (\ref{Eq: DCS4}), the combination of (\ref{Eq: new 1}) and (\ref{Eq: new 2}), and the combination of (\ref{Eq: DCS1}) and (\ref{Eq: DCS2}).
The latter three mixed combinations are referred to as the \emph{mixed discrete conformal structure I}, \emph{mixed discrete conformal structure II}, and \emph{mixed discrete conformal structure III}, respectively.
\end{theorem}

\begin{remark}\label{Rmk: 1}
As noted in \cite[Remark 1.5]{X-Z DCS1}, for simplicity, we may always assume that $\alpha: B\rightarrow \{-1,0,1\}$ in the discrete conformal structures (\ref{Eq: DCS3}) and (\ref{Eq: DCS4}), $\alpha\equiv -1$ in the discrete conformal structures (\ref{Eq: new 1}) and (\ref{Eq: new 2}), and $C\equiv0$ in the discrete conformal structures (\ref{Eq: DCS1}) and (\ref{Eq: DCS2}).
\end{remark}

The main result in \cite{X-Z DCS2} is the following theorem, which establishes the global rigidity and existence of discrete conformal structures on surfaces with boundary.

\begin{theorem}[\cite{X-Z DCS2}, Theorem 1.8]\label{Thm: rigidity and image}
Let $(\Sigma,\mathcal{T})$ be an ideally triangulated surface with boundary, and let $d=d(f)$ be a discrete conformal structure on $(\Sigma,\mathcal{T})$.
\begin{description}
\item[(i)]
For the discrete conformal structure (\ref{Eq: DCS3}), let $\alpha: B\rightarrow \{-1,0,1\}$ and $\eta\in \mathbb{R}_{>0}^E$ be the weights on $(\Sigma,\mathcal{T})$ satisfying $\eta_{ij}>\alpha_i\alpha_j$ for any two adjacent boundary components $i,j\in B$.
Then the discrete conformal factor $f$ is uniquely determined by its generalized combinatorial curvature $K\in \mathbb{R}^N_{>0}$.
Furthermore, if $\alpha: B\rightarrow \{0,1\}$, then
the image of $K$ is $\mathbb{R}^N_{>0}$.

\item[(ii)]
For the discrete conformal structure (\ref{Eq: new 1}), let $\eta\in (-1,+\infty)^E$ be the weight on $(\Sigma,\mathcal{T})$.
Then the discrete conformal factor $f$ is uniquely determined by its generalized combinatorial curvature $K\in \mathbb{R}^N_{>0}$.
Furthermore, if $\eta\in (-1,0]^E$, then
the image of $K$ is $\mathbb{R}^N_{>0}$.

\item[(iii)]
For the discrete conformal structure (\ref{Eq: DCS1}), let $\eta\in \mathbb{R}_{>0}^E$ be the weight on $(\Sigma,\mathcal{T})$.
Then the discrete conformal factor $f$ is uniquely determined by its generalized combinatorial curvature $K\in \mathbb{R}^N_{>0}$.
Furthermore, the image of $K$ is $\mathbb{R}^N_{>0}$.

\item[(iv)]
For the mixed discrete conformal structure I,
let $\alpha: B\rightarrow \{-1,0,1\}$ and $\eta\in (1,+\infty)^E$ be the weights on $(\Sigma,\mathcal{T})$.
Then the discrete conformal factor $f$ is uniquely determined by its generalized combinatorial curvature $K\in \mathbb{R}^N_{>0}$.
Furthermore, if $\alpha: B\rightarrow \{0,1\}$, and $\alpha_j$ and $\alpha_k$ cannot both be 1 simultaneously for any right-angled hyperbolic hexagon with edge lengths $l_{ij},l_{ik}$ given by (\ref{Eq: DCS4}) and $l_{jk}$ given by (\ref{Eq: DCS3}),
then the image of $K$ is $\mathbb{R}^N_{>0}$.

\item[(v)]
For the mixed discrete conformal structure II, let $\eta\in [1,+\infty)^E$ be the weight on $(\Sigma,\mathcal{T})$.
Then the discrete conformal factor $f$ is uniquely determined by its generalized combinatorial curvature $K\in \mathbb{R}^N_{>0}$.

\item[(vi)]
For the mixed discrete conformal structure III, let $\eta\in \mathbb{R}_{>0}^E$ be the weight on $(\Sigma,\mathcal{T})$.
Then the discrete conformal factor $f$ is uniquely determined by its generalized combinatorial curvature $K\in \mathbb{R}^N_{>0}$.
Furthermore, the image of $K$ is $\mathbb{R}^N_{>0}$.
\end{description}
\end{theorem}

The proof of the rigidity part in Theorem \ref{Thm: rigidity and image} consists of three steps, of which the first two are as follows.
The first step is to characterize the admissible space of discrete conformal factors $f$ for a right-angled hyperbolic hexagon.
Through the variable transformation $u_i=u_i(f_i)$, the admissible space $\Omega_{ijk}$ of $f$ is transformed into the admissible space $\mathcal{U}_{ijk}$ of $u$, which is proven to be convex and simply connected.
The second step involves proving that the Jacobian of the generalized angles $\theta$ with respect to $u$ for a right-angled hyperbolic hexagon is symmetry and negative definiteness.
As a direct consequence, the Jacobian of the generalized combinatorial curvature $K$ with respect to $u$ is also symmetric and negative definite on the admissible space $\mathcal{U}(\eta) = \bigcap_{\{ijk\} \in F}\, \mathcal{U}_{ijk}(\eta)$.

A critical component of the proof hinges on constructing the necessary variable transformation $u_i=u_i(f_i)$.
To establish this transformation, we first refer to the variation formula of the generalized angles $\theta_i^{jk}$ with respect to $f_j$ for a right-angled hyperbolic hexagon $\{ijk\}\in F$, as derived in Lemma 2.17 of \cite{X-Z DCS2}.
Notably, the asymmetry of this formula in the indices $i$ and $j$ motivates the introduction of the variable transformation $u_i=u_i(f_i)$, the explicit expressions and associated discussions of which are provided in Lemma 2.21 and Remark 2.22 of \cite{X-Z DCS2}.

\subsection{Combinatorial curvature flows on surfaces with boundary}

\begin{definition}\label{Def: CCF}
Under the same assumptions as those in Theorem \ref{Thm: rigidity and image},
let $\overline{K}\in \mathbb{R}_{>0}^N$ be a given function defined on $B=\{1,2,...,N\}$.
For the discrete conformal structures on an ideally triangulated surface with boundary $(\Sigma,\mathcal{T})$,
the combinatorial Ricci flow is defined as
\begin{eqnarray}\label{Eq: CRF}
\begin{cases}
\frac{du_i}{dt}=K_i-\overline{K}_i,\\
u_i(0)=u_0,
\end{cases}
\end{eqnarray}
and the combinatorial Calabi flow is defined as
\begin{eqnarray}\label{Eq: CCF}
\begin{cases}
\frac{du_i}{dt}=-\Delta(K-\overline{K})_i,\\
u_i(0)=u_0,
\end{cases}
\end{eqnarray}
where $\Delta=(\frac{\partial K_i}{\partial u_j})_{N\times N}$ is the discrete Laplace operator.
\end{definition}

\begin{remark}
The combinatorial Ricci flow and the combinatorial Calabi flow introduced in Definition \ref{Def: CCF} unify and generalize various known forms of combinatorial Ricci/Yamabe flows and combinatorial Calabi flows for specific types of discrete conformal structures on ideally triangulated surfaces with boundary.
More precisely, for the discrete conformal structure (\ref{Eq: DCS3}),
when $\alpha\equiv1$, the combinatorial Ricci flow (\ref{Eq: CRF}) reduces to the generalized combinatorial curvature flow introduced by Guo-Luo in \cite{GL2},while the combinatorial Calabi flow (\ref{Eq: CCF}) reduces to the combinatorial Calabi flow we proposed for generalized circle packings of $(-1,-1,-1)$-type in \cite{XZ};
when $\alpha\equiv0$, the combinatorial Ricci flow (\ref{Eq: CRF}) reduces to both Guo's combinatorial Yamabe flow for vertex scaling on surfaces with boundary in \cite{Guo}, and the modified combinatorial Yamabe flow for Guo's vertex scaling on surfaces with boundary introduced by Li-Xu-Zhou in \cite{Li-Xu-Zhou}, while the combinatorial Calabi flow (\ref{Eq: CCF}) reduces to the combinatorial Calabi flow for Guo's vertex scaling on surfaces with boundary introduced by Luo-Xu in \cite{Luo-Xu};
when $\alpha\equiv-1$, both the combinatorial Ricci flow (\ref{Eq: CRF}) and the combinatorial Calabi flow (\ref{Eq: CCF}) reduce to their counterparts in Xu's discrete conformal structures on surfaces with boundary in \cite{Xu22}.
\end{remark}

The main result of this paper is the following theorem, which establishes the longtime existence and convergence of solutions to the combinatorial Ricci flow (\ref{Eq: CRF}) and the combinatorial Calabi flow (\ref{Eq: CCF}) for various discrete conformal structures on surfaces with boundary.

\begin{theorem}\label{Thm: converge}
Let $(\Sigma,\mathcal{T})$ be an ideally triangulated surface with boundary, and let $d=d(f)$ be a discrete conformal structure on $(\Sigma,\mathcal{T})$.
\begin{description}
\item[(i)]
For the discrete conformal structure (\ref{Eq: DCS3}), let $\alpha: B\rightarrow \{0,1\}$ and $\eta\in \mathbb{R}_{>0}^E$ be the weights on $(\Sigma,\mathcal{T})$ satisfying $\eta_{ij}>\alpha_i\alpha_j$ for any two adjacent boundary components $i,j\in B$.
Then the solutions to the combinatorial Ricci flow (\ref{Eq: CRF}) and the combinatorial Calabi flow (\ref{Eq: CCF}) exist for all time and converge exponentially fast for any initial value $u(0)$, respectively.

\item[(ii)]
For the discrete conformal structure (\ref{Eq: DCS1}), let $\eta\in \mathbb{R}_{>0}^E$ be the weight on $(\Sigma,\mathcal{T})$.
Then the solutions to the combinatorial Ricci flow (\ref{Eq: CRF}) and the combinatorial Calabi flow (\ref{Eq: CCF}) exist for all time and converge exponentially fast for any initial value $u(0)$, respectively.

\item[(iii)]
For the discrete conformal structure (\ref{Eq: new 1}), and the mixed discrete conformal structures I and III, under the assumptions of Theorem \ref{Thm: rigidity and image} (ii), (iv) and (vi),
there exists a constant $\delta>0$ such that if $\Vert K(u(0))-\overline{K}\Vert<\delta$,
then the solutions to the combinatorial Ricci flow (\ref{Eq: CRF}) and the combinatorial Calabi flow (\ref{Eq: CCF})  exist for all time and converge exponentially fast, respectively.

\item[(iv)]
For the mixed discrete conformal structure II, let $\eta\in [1,+\infty)^E$ be the weight on $(\Sigma,\mathcal{T})$.
If the solution $u(t)$ to the combinatorial Ricci flow (\ref{Eq: CRF}) or the combinatorial Calabi flow (\ref{Eq: CCF}) converges to $\overline{u}\in \mathcal{U}(\eta)$, then $K(\overline{u})=\overline{K}$.
Conversely, suppose that there exists a discrete conformal factor $\overline{u}$ on $(\Sigma,\mathcal{T})$ satisfying $K(\overline{u})=\overline{K}$, then there exists a constant $\delta>0$ such that if $\Vert K(u(0))-\overline{K}\Vert<\delta$, the solutions to the combinatorial Ricci flow (\ref{Eq: CRF}) and the combinatorial Calabi flow (\ref{Eq: CCF})  exist for all time and converge exponentially fast to $\overline{u}$.
\end{description}
\end{theorem}

\begin{remark}
If $\alpha\equiv1$, then Theorem \ref{Thm: converge} (i) generalizes our previous results in \cite{XZ}.
If $\alpha\equiv0$, then Theorem \ref{Thm: converge} (i) generalizes Guo's results in \cite{Guo} and Li-Xu-Zhou's results in \cite{Li-Xu-Zhou} for the combinatorial Ricci flow (\ref{Eq: CRF}), and Luo-Xu's results in \cite{Luo-Xu} for the combinatorial Calabi flow (\ref{Eq: CCF}).
\end{remark}

\subsection{Organization of the paper}
In Section \ref{Sec: CCF}, we investigate fundamental properties of the combinatorial curvature flows, and introduce two significant functions.
In Section \ref{Sec: DSC3}, we analyze the combinatorial curvature flows for the discrete conformal structure (\ref{Eq: DCS3}) and prove Theorem \ref{Thm: converge} (i).
In Section \ref{Sec: DSC1}, we consider the combinatorial curvature flows of the discrete conformal structure (\ref{Eq: DCS1}) and prove Theorem \ref{Thm: converge} (ii).
In Section \ref{Sec: remaining}, we consider the combinatorial curvature flows for the remaining discrete conformal structures and prove Theorem \ref{Thm: converge} (iii) and (iv).
\\
\\
\textbf{Acknowledgements}\\[8pt]
The research of X. Xu is supported by National Natural Science Foundation of China under grant no. 12471057.

\section{Combinatorial curvature flows on surfaces with boundary}\label{Sec: CCF}

Since the combinatorial Ricci flow (\ref{Eq: CRF}) and the combinatorial Calabi flow (\ref{Eq: CCF}) are ODE systems with smooth coefficients,
standard ODE theory ensures the local existence of their solutions around the initial time $t=0$.
By applying the Lyapunov Stability Theorem (see \cite[Chapter 5]{Pontryagin}),
we can establish the longtime existence and convergence of solutions to both the combinatorial Ricci flow (\ref{Eq: CRF}) and the combinatorial Calabi flow (\ref{Eq: CCF}) for initial data with small energy.

For general initial data, to establish the longtime existence and convergence of solutions to the combinatorial Ricci flow (\ref{Eq: CRF}) and the combinatorial Calabi flow (\ref{Eq: CCF}),
it suffices to show that these solutions stay in a compact subset of the admissible space, i.e., the solutions cannot reach the boundary of the admissible space.
Consequently, we need to analyze the behavior of these combinatorial curvature flows as they approach the boundary of the admissible space.

We first introduce the following two useful energy functions:
\begin{equation}\label{Eq: EF1}
\mathcal{E}(u)
=-\int_0^{u}\sum_{i=1}^{N}(K_i-\overline{K}_i)du_i,
\end{equation}
and
\begin{equation}\label{Eq: EF2}
\mathcal{C}(u)
=\frac{1}{2}\Vert K-\overline{K}\Vert^2
=\frac{1}{2}\sum_{i=1}^{N}(K_i-\overline{K}_i)^2.
\end{equation}

Note that the function $\mathcal{E}(u)$ is a well-defined strictly convex function on the admissible space $\mathcal{U}(\eta)$ of $u$.
This is due to the fact that, for different discrete conformal structures, the admissible spaces of $u$ are convex polytopes, and the Jacobians $\frac{\partial (K_i,..., K_N)}{\partial(u_i,...,u_N)}$ are symmetric and negative definite.
Specifically, the admissible spaces corresponding to the discrete conformal structures (\ref{Eq: DCS3}) and (\ref{Eq: DCS1}) are presented in Theorem \ref{Thm: ASC 1} and Theorem \ref{Thm: ASC3}, respectively.
Moreover, the Jacobians $\frac{\partial (K_i,..., K_N)}{\partial(u_i,...,u_N)}$ for these two discrete conformal structures are derived in Theorem \ref{Thm: matrix negative 1} and Theorem \ref{Thm: matrix negative 2}, respectively.
For the remaining discrete conformal structures, the admissible spaces are detailed in Section 4, Corollary 6.2, Corollary 7.2, and Corollary 8.2 of \cite{X-Z DCS2}, with the associated Jacobians given in Theorem 4.1, Corollary 6.4, Corollary 7.4, and Corollary 8.4 of \cite{X-Z DCS2}.

\begin{proposition}\label{Prop: property}
The functions $\mathcal{E}(u)$ and $\mathcal{C}(u)$ have the following properties.
\begin{description}
\item[(i)]
The combinatorial Ricci flow (\ref{Eq: CRF}) is the negative gradient flow of $\mathcal{E}(u)$,
and the combinatorial Calabi flow (\ref{Eq: CCF}) is the negative gradient flow of $\mathcal{C}(u)$.
\item[(ii)]
Both functions $\mathcal{E}(u)$ and $\mathcal{C}(u)$ are decreasing along both the combinatorial Ricci flow (\ref{Eq: CRF}) and the combinatorial Calabi flow (\ref{Eq: CCF}).
\end{description}
\end{proposition}
\proof
For the part (1), the conclusion is derived from the following calculations:
\begin{gather*}
\nabla_{u_i}\mathcal{E}(u)=-(K_i-\overline{K}_i)
=-\frac{du_i}{dt},\\
\nabla_{u_i}\mathcal{C}(u)=\sum_{i=1}^N\frac{\partial K_j}{\partial u_i}(K_j-\overline{K}_j)
=\Delta(K-\overline{K})_i=-\frac{du_i}{dt}.
\end{gather*}

For the part (2), a routine calculation yields
\begin{equation*}
\frac{d\mathcal{E}(u(t))}{dt}
=\sum_{i=1}^{N}\frac{\partial \mathcal{E}}{\partial u_i}\frac{d u_i}{dt}
=-\sum_{i=1}^{N}(K_i-\overline{K}_i)^2\leq0.
\end{equation*}
This implies that $\mathcal{E}(u)$ is decreasing along the combinatorial Ricci flow (\ref{Eq: CRF}). Similarly,
\begin{equation*}
\frac{d\mathcal{C}(u(t))}{dt}
=\sum_{i,j=1}^{N}\frac{\partial \mathcal{C}}{\partial K_i}\frac{\partial K_i}{\partial u_j}\frac{d u_j}{dt}
=\sum_{i,j=1}^{N}(K_i-\overline{K}_i)\frac{\partial K_i}{\partial u_j}(K_j-\overline{K}_j)
=(K-\overline{K})^T\Delta(K-\overline{K})\leq0,
\end{equation*}
where the inequality follows from the negative definiteness of $\Delta=(\frac{\partial K_i}{\partial u_j})_{N\times N}$, and the right-hand side is strictly negative unless $K=\overline{K}$.
This implies that $\mathcal{C}(u)$ is decreasing along the combinatorial Ricci flow (\ref{Eq: CRF}).
Similar calculations yield
\begin{gather*}
\frac{d\mathcal{E}(u(t))}{dt}
=\sum_{i=1}^{N}\frac{\partial \mathcal{E}}{\partial u_i}\frac{d u_i}{dt}
=\sum_{i=1}^{N}(K-\overline{K})_i\Delta(K-\overline{K})_i
\leq0,\\
\frac{d\mathcal{C}(u(t))}{dt}
=\sum_{i=1}^{N}\frac{\partial \mathcal{C}}{\partial u_i}\frac{d u_i}{dt}
=-\sum_{i=1}^{N}(\Delta(K-\overline{K})_i)^2\leq0.
\end{gather*}
This completes the proof.
\qed

\section{Combinatorial curvature flows for the discrete conformal structure (\ref{Eq: DCS3})}\label{Sec: DSC3}

Suppose $(\Sigma,\mathcal{T},\alpha,\eta)$ is a weighted triangulated surface with boundary,
where the weights are given by $\alpha: B\rightarrow \{-1,0,1\}$ and $\eta\in \mathbb{R}^E_{>0}$.
By Lemma 2.21 and Remark 2.22 in \cite{X-Z DCS2}, for any $i\in B$,
\begin{eqnarray}\label{Eq: F1}
u_i=
\begin{cases}
f_i,  &{\alpha_i=0},\\
\frac{1}{2}\ln\left|\frac{\sqrt{1+\alpha_ie^{2f_i}}-1}
{\sqrt{1+\alpha_ie^{2f_i}}+1}\right|,
&{\alpha_i=\pm1}.
\end{cases}
\end{eqnarray}
For a right-angled hyperbolic hexagon $\{ijk\}\in F$ with edge lengths $l_{ij}, l_{jk}, l_{ki}$ given by (\ref{Eq: DCS3}),
the admissible space of the discrete conformal factors $u$ is defined as
\begin{equation*}
\mathcal{U}_{ijk}(\eta)
=\{(u_i,u_j,u_k)\in \mathbb{R}^{N_1}\times \mathbb{R}_{<0}^{3-N_1} \mid l_{ij}>0, l_{jk}>0, l_{ki}>0\},
\end{equation*}
where $N_1$ is the number of the boundary components $i\in B$ with $\alpha_i=0$.
Furthermore, we have the following results.

\begin{theorem}[\cite{X-Z DCS2}, Theorem 3.1]
\label{Thm: ASC 1}
Suppose $(\Sigma,\mathcal{T},\alpha,\eta)$ is a weighted triangulated surface with boundary, where the weighs are given by $\alpha: B\rightarrow \{-1,0,1\}$ and $\eta\in \mathbb{R}_{>0}^E$ satisfying $\eta_{ij}>\alpha_i\alpha_j$ for any two adjacent boundary components $i,j\in B$.
Then the admissible space
\begin{equation}\label{Eq: admissible space a2}
\mathcal{U}_{ijk}(\eta)
=\{(u_i,u_j,u_k)\in \mathbb{R}^{N_1}\times \mathbb{R}_{<0}^{3-N_1}\mid u_r+u_s>C(\eta_{rs}),\ \forall\{r,s\}\subseteq\{i,j,k\}\}
\end{equation}
is a convex polytope,
where $C(\eta_{rs})$ is a constant depending on $\eta_{rs}$.
As a result, the admissible space
$\mathcal{U}(\eta)=\bigcap_{\{ijk\}\in F}\mathcal{U}_{ijk}(\eta)$
is also a convex polytope.
\end{theorem}

\begin{theorem}[\cite{X-Z DCS2}, Theorem 3.2 and Corollary 3.3]\label{Thm: matrix negative 1}
Under the same assumptions as those in Theorem \ref{Thm: ASC 1},
for a right-angled hyperbolic hexagon $\{ijk\}\in F$ on $\mathcal{U}_{ijk}(\eta)$,
the Jacobian  $\frac{\partial(\theta^{jk}_i,\theta^{ik}_j,\theta^{ij}_k)}
{\partial(u_i,u_j,u_k)}$ is symmetric and negative definite.
As a result, the Jacobian $\frac{\partial (K_i,..., K_N)}{\partial(u_i,...,u_N)}$ is symmetric and negative definite on $\mathcal{U}(\eta)$.
\end{theorem}

The combinatorial curvature flows for the discrete conformal structure (\ref{Eq: DCS3}) is redefined as follows.

\begin{definition}\label{Def: CCF1}
Suppose $(\Sigma,\mathcal{T},\alpha,\eta)$ is a weighted triangulated surface with boundary, where the weighs are given by $\alpha: B\rightarrow \{0,1\}$ and $\eta\in \mathbb{R}_{>0}^E$ satisfying $\eta_{ij}>\alpha_i\alpha_j$ for any two adjacent boundary components $i,j\in B$.
Let $\overline{K}\in \mathbb{R}_{>0}^N$ be a given function defined on $B=\{1,2,...,N\}$.
The combinatorial Ricci flow for the discrete conformal structure (\ref{Eq: DCS3}) is defined as
\begin{eqnarray}\label{Eq: CRF1}
\begin{cases}
\frac{du_i}{dt}=K_i-\overline{K}_i,\\
u_i(0)=u_0.
\end{cases}
\end{eqnarray}
The combinatorial Calabi flow for the discrete conformal structure (\ref{Eq: DCS3}) is defined as
\begin{eqnarray}\label{Eq: CCF1}
\begin{cases}
\frac{du_i}{dt}=-\Delta(K-\overline{K})_i,\\
u_i(0)=u_0,
\end{cases}
\end{eqnarray}
where $\Delta=(\frac{\partial K_i}{\partial u_j})_{N\times N}$ is the discrete Laplace operator in Theorem \ref{Thm: matrix negative 1}.
\end{definition}

By (\ref{Eq: admissible space a2}), let $\alpha: B\rightarrow \{0,1\}$ and let $N_1$ is the number of the boundary components $i\in B$ with $\alpha_i=0$.
The boundary of the admissible space $\mathcal{U}(\eta)$ in $[-\infty,+\infty]^{N_1}\times[-\infty,0]^{N-N_1}$ comprises three distinct parts:
\begin{description}
\item[(i)]
There exists $i\in B$ such that either $u_i=\pm\infty$ with $\alpha_i=0$ or $u_i=-\infty$ with $\alpha_i=1$, the set of which is denoted by $\partial_{\infty}\mathcal{U}(\eta)$;
\item[(ii)]
There exists $i\in B$ such that $u_i=0$ with $\alpha_i=1$, the set of which is denoted by $\partial_{0}\mathcal{U}(\eta)$;
\item[(iii)]
There exists an ideal edge $\{ij\}\in E$ such that $l_{ij}=0$, which is equivalent to $u_i+u_j=C(\eta_{ij})$ for the admissible space $\mathcal{U}_{ijk}(\eta)$, the set of which is denoted by $\partial_{ij}\mathcal{U}(\eta)$.
\end{description}

\subsection{Combinatorial Ricci flow for the discrete conformal structure (\ref{Eq: DCS3})}

The following three lemmas show that the solution $u(t)$ to the combinatorial Ricci flow (\ref{Eq: CRF1}) cannot reach the boundaries $\partial_{\infty}\mathcal{U}(\eta)$, $\partial_{0}\mathcal{U}(\eta)$, and $\partial_{l}\mathcal{U}(\eta)=\bigcup_{\{ij\}\in E}\partial_{ij}\mathcal{U}(\eta)$, respectively.

\begin{lemma}\label{Lem: boundary a1}
Under the same assumptions as those in Definition \ref{Def: CCF1},
the solution $u(t)$ to the combinatorial Ricci flow (\ref{Eq: CRF1}) stays in a bounded subset of $\mathbb{R}^{N_1}\times \mathbb{R}_{<0}^{N-N_1}$.
Consequently, it cannot reach the boundary $\partial_{\infty}\mathcal{U}(\eta)$.
\end{lemma}
\proof
By Theorem \ref{Thm: rigidity and image} (i), there exists $\overline{u}\in \mathcal{U}(\eta)$ such that $K(\overline{u})=\overline{K}$.
It follows that $\nabla\mathcal{E}(\overline{u})
=-(K-\overline{K})|_{u=\overline{u}}=0$.
By Theorem \ref{Thm: matrix negative 1},
$\mathrm{Hess}_u\, \mathcal{E} =-\Delta>0$,
which implies $\mathcal{E}(u)$ is a strictly convex function on $\mathcal{U}(\eta)\subseteq\mathbb{R}^{N_1}\times \mathbb{R}_{<0}^{N-N_1}$.
This strict convexity ensures that
$\lim_{\Vert u\Vert\rightarrow +\infty}\mathcal{E}(u)=+\infty$,
so $\mathcal{E}(u)$ is a proper function.
Moreover, Proposition \ref{Prop: property} (ii) shows that $\mathcal{E}(u)$ is decreasing along the combinatorial Ricci flow (\ref{Eq: CRF1}).
Consequently, $\mathcal{E}(u(t))\leq \mathcal{E}(u(0))$, meaning that the solution $u(t)$ to the combinatorial Ricci flow \eqref{Eq: CRF1} stays in a bounded subset of $\mathbb{R}^{N_1}\times \mathbb{R}_{<0}^{N-N_1}$.
\qed

\begin{lemma}\label{Lem: boundary a2}
Under the same assumptions as those in Definition \ref{Def: CCF1},
if $\alpha_i=1$, then the solution $u_i(t)$ to the combinatorial Ricci flow (\ref{Eq: CRF1}) is uniformly bounded from above in $\mathbb{R}_{<0}$.
Consequently, it cannot reach the boundary $\partial_{0}\mathcal{U}(\eta)$.
\end{lemma}
\proof
By (\ref{Eq: F1}), if $\alpha_i=1$, then
$-\coth u_i=\sqrt{1+e^{2f_i}}.$
Hence, as $\lim_{t\rightarrow T}u_i(t)=0^-$ for some $T\in(0,+\infty]$,
we have $\lim_{t\rightarrow T}f_i(t)=+\infty$.
By Lemma 3.6 of \cite{X-Z DCS2}, $\theta^{jk}_{i}\rightarrow 0^-$ uniformly as $f_i\rightarrow+\infty$.
Consequently, there exists a constant $c\in (-\infty,0)$ such that whenever $u_i(t)>c$, the hyperbolic arc length $\theta^{jk}_{i}$ is less than any given $\epsilon>0$,
which implies $K_i<\overline{K}_i$.
Since $\lim_{t\rightarrow T}u_i(t)=0^-$, we may choose a time $t_0\in (0,T)$ such that $u_i(t_0)>c$.
Set $a=\inf\{t<t_0\mid u_i(s)>c, \forall s\in (t,t_0]\}$, and then $u_i(a)=c$.
Note that $\frac{du_i}{dt}=K_i-\overline{K}_i<0$ on $(a,t_0]$, it follows that $u_i(t_0)<u_i(a)=c$, which contradicts $u_i(t_0)>c$.
Therefore, $u_i(t)$ is uniformly bounded from above in $\mathbb{R}_{<0}$.
\qed

\begin{lemma}\label{Lem: boundary a3}
Under the same assumptions as those in Definition \ref{Def: CCF1},
the solution $u(t)$ to the combinatorial Ricci flow (\ref{Eq: CRF1}) cannot reach the boundary $\partial_{l}\mathcal{U}(\eta)$.
\end{lemma}
\proof
If there exists an ideal edge $\{ij\}\in E$ such that $l_{ij}=0$, i.e., $u_i+u_j=C(\eta_{ij})$, then
\begin{equation*}
\cosh \theta^{ik}_j=\frac{\cosh l_{ik}+\cosh l_{ij}\cosh l_{jk}}{\sinh l_{ij}\sinh l_{jk}}>\frac{\cosh l_{ij}\cosh l_{jk}}{\sinh l_{ij}\sinh l_{jk}}\geq\frac{\cosh l_{ij}}{\sinh l_{ij}}\rightarrow +\infty.
\end{equation*}
By the definition of $K$ in (\ref{Eq: K}), i.e., $K_j=\sum_{\{ijk\}\in F}\theta^{ik}_j$,
it follows that $K_j\rightarrow +\infty$.
This implies $\mathcal{C}(u)\rightarrow +\infty$ by (\ref{Eq: EF2}).
However, Proposition \ref{Prop: property} (ii) shows that $\mathcal{C}(u)$ is decreasing along the combinatorial Ricci flow (\ref{Eq: CRF1}),
which implies $\mathcal{C}(u(t))\leq \mathcal{C}(u(0))$.
This leads to a contradiction.
Therefore, the solution $u(t)$ to the combinatorial Ricci flow (\ref{Eq: CRF1}) cannot reach the boundary $\partial_{l}\mathcal{U}(\eta)$.
\qed

As an immediate corollary of Lemma \ref{Lem: boundary a1}, Lemma \ref{Lem: boundary a2} and Lemma \ref{Lem: boundary a3},
we conclude that the solution $u(t)$ to the combinatorial Ricci flow \eqref{Eq: CRF1} exists for all time,
which corresponds to the existence part of Theorem \ref{Thm: converge} (i).

\begin{corollary}\label{Cor: exist a1}
Under the same assumptions as those in Definition \ref{Def: CCF1}, the solution $u(t)$ to the combinatorial Ricci flow \eqref{Eq: CRF1} stays in a compact subset of $\mathcal{U}(\eta)$.
As a result, this solution exists for all time.
\end{corollary}

The following theorem establishes the convergence of the solution $u(t)$ to the combinatorial Ricci flow \eqref{Eq: CRF1},
which corresponds to the convergence part of Theorem \ref{Thm: converge} (i).

\begin{theorem}\label{Thm: converge a1}
Under the same assumptions as those in Definition \ref{Def: CCF1},
the solution $u(t)$ to the combinatorial Ricci flow \eqref{Eq: CRF1} converges exponentially fast.
\end{theorem}
\proof
By Theorem \ref{Thm: rigidity and image} (i), there exists $\overline{u}\in \mathcal{U}(\eta)$ such that $K(\overline{u})=\overline{K}$.
Combining the fact that $\mathcal{E}(u(t))$ is decreasing along the combinatorial Ricci flow (\ref{Eq: CRF1}) in Proposition \ref{Prop: property} (ii), and the result that the solution $u(t)$ to the combinatorial Ricci flow (\ref{Eq: CRF1}) stays in a compact subset of $\mathcal{U}(\eta)$ in Corollary \ref{Cor: exist a1},
we conclude that $\lim_{t\rightarrow +\infty}\mathcal{E}(u(t))$ exists.
Consequently, there exists a sequence $\xi_n\in(n,n+1)$ such that as $n\rightarrow +\infty$,
\begin{equation*}
\mathcal{E}(u(n+1))-\mathcal{E}(u(n))
=(\mathcal{E}(u(t))'|_{\xi_n}=\nabla \mathcal{E}\cdot\frac{du_i}{dt}|_{\xi_n}
=-\sum_{i=1}^{N}(K_i(u(\xi_n))-\overline{K}_i)^2
\rightarrow 0.
\end{equation*}
Hence, $\lim_{n\rightarrow +\infty}K_i(u(\xi_n))=\overline{K}_i=K_i(\overline{u})$ for all $i\in B$.
Since $\{u(t)\}\subset\subset \mathcal{U}(\eta)$ by Corollary \ref{Cor: exist a1},
there exists $u^*\in \mathcal{U}(\eta)$ and a subsequence of $\{u(\xi_n)\}$, still denoted by $\{u(\xi_n)\}$ for simplicity,
such that $\lim_{n\rightarrow \infty}u(\xi_n)=u^*$.
This implies $K_i(u^*)=\lim_{n\rightarrow +\infty}K_i(u(\xi_n))=K_i(\overline{u})$.
It follows from \ref{Thm: rigidity and image} (i) that $u^*=\overline{u}$.
Therefore, $\lim_{n\rightarrow \infty}u(\xi_n)=\overline{u}$.

Set $\Gamma(u)=K-\overline{K}$.
By Theorem \ref{Thm: matrix negative 1}, the matrix $D\Gamma|_{u=\overline{u}}=\Delta<0$, i.e., it has $N$ negative eigenvalues.
This implies that $\overline{u}$ is a local attractor of the combinatorial Ricci flow (\ref{Eq: CRF1}).
The conclusion follows from Lyapunov Stability Theorem (\cite[Chapter 5]{Pontryagin}).
\qed

\subsection{Combinatorial Calabi flow for the discrete conformal structure (\ref{Eq: DCS3})}

The following lemma shows that the solution $u(t)$ to the combinatorial Calabi flow (\ref{Eq: CCF1}) cannot reach the boundaries $\partial_{\infty}\mathcal{U}(\eta)$ and $\partial_{l}\mathcal{U}(\eta)$.
Since the proof is nearly identical to those of Lemma \ref{Lem: boundary a1} and Lemma \ref{Lem: boundary a3}, it is omitted here.

\begin{lemma}\label{Lem: boundary a4}
Under the same assumptions as those in Definition \ref{Def: CCF1}, the solution $u(t)$ to the combinatorial Calabi flow (\ref{Eq: CCF1}) stays in a bounded subset of $\mathbb{R}^{N_1}\times \mathbb{R}_{<0}^{N-N_1}$ and cannot reach the boundary $\partial_{l}\mathcal{U}(\eta)$.
\end{lemma}

To prove the solution $u(t)$ to the combinatorial Calabi flow (\ref{Eq: CCF1}) cannot reach the boundary $\partial_{0}\mathcal{U}(\eta)$,
we need the following lemma.

\begin{lemma}\label{lem: CCF a1}
Under the same assumptions as those in Definition \ref{Def: CCF1},
let $\{ijk\}\in F$ be a hyperbolic right-angled hexagon with $\alpha_i=1$.
For any constant $C\in \mathbb{R}$, there exists a constant $W=W(C)>0$ such that if $f_i\geq W$, then
\begin{equation}\label{Eq: F53}
\bigg|\frac{\partial \theta^{jk}_i}{\partial u_i}\bigg|
>C\bigg(\bigg|\frac{\partial \theta^{jk}_i}{\partial u_j}\bigg|+\bigg|\frac{\partial \theta^{jk}_i}{\partial u_k}\bigg|\bigg).
\end{equation}
\end{lemma}
\proof
For simplicity, set $l_r=l_{st}$, $\theta_r=\theta^{st}_r$, $A=\sinh l_r\sinh l_s\sinh \theta_t$, $S_r=e^{f_r}$, and $C_r=\sqrt{1+\alpha_re^{2f_r}}$,
where $\{r,s,t\}=\{i,j,k\}$.
A routine calculation yields
\begin{equation}\label{Eq: F4}
\begin{aligned}
\frac{\partial(\theta_i,\theta_j,\theta_k)}{\partial(u_i,u_j,u_k)}
&=\frac{\partial(\theta_i,\theta_j,\theta_k)}{\partial(l_i,l_j,l_k)}\cdot
\frac{\partial(l_i,l_j,l_k)}{\partial(f_i,f_j,f_k)}\cdot
\frac{\partial(f_i,f_j,f_k)}{\partial(u_i,u_j,u_k)}\\
&=-\frac{1}{A}
 \left(
   \begin{array}{ccc}
     \sinh l_i & 0 & 0 \\
     0 & \sinh l_j & 0 \\
     0 & 0 & \sinh l_k \\
   \end{array}
 \right)
 \left(
   \begin{array}{ccc}
     -1 & \cosh\theta_k & \cosh\theta_j \\
     \cosh\theta_k & -1 & \cosh\theta_i \\
     \cosh\theta_j & \cosh\theta_i & -1 \\
   \end{array}
 \right)\\
& \ \ \ \ \ \ \ \times
 \left(
   \begin{array}{ccc}
     0 & \coth d_{jk} & \coth d_{kj} \\
     \coth d_{ik} & 0 & \coth d_{ki} \\
     \coth d_{ij} & \coth d_{ji} & 0 \\
   \end{array}
 \right)
 \left(
   \begin{array}{ccc}
     C_i & 0 & 0 \\
     0 & C_j & 0 \\
     0 & 0 & C_k \\
   \end{array}
 \right).
\end{aligned}
\end{equation}
One can refer to Subsection 3.2 in \cite{X-Z DCS2} for (\ref{Eq: F4}).
Moreover, the formula (\ref{Eq: DCS3}) can be written as $\cosh l_k=-C_iC_j+\eta_{ij}S_iS_j$,
and the formula (\ref{Eq: d3}) can be written as
\begin{equation}\label{Eq: F5}
\coth d_{ij}
=\frac{-\alpha_iS^2_i C_j}{\sinh l_k C_i}
+\frac{\eta_{ij}S_iS_j}{\sinh l_k}
=\frac{-\alpha_iS^2_i C_j+C_i(\cosh l_k+C_iC_j)}{\sinh l_k C_i}
=\frac{C_j+\cosh l_k C_i}{\sinh l_k C_i},
\end{equation}
where $C_i^2-\alpha_iS^2_i=1$ is used in the final equality.
Substituting (\ref{Eq: F5}) together with $\cosh \theta_i=\frac{\cosh l_i+\cosh l_j\cosh l_k}{\sinh l_j\sinh l_k}$ into (\ref{Eq: F4}) yields
\begin{equation}\label{Eq: F54}
\begin{aligned}
\frac{\partial \theta_i}{\partial u_i}
=-\frac{1}{A}[&\frac{1}{\sinh^2 l_j}(C_k+C_i\cosh l_j)(\cosh l_k+\cosh l_i\cosh l_j)\\
&+\frac{1}{\sinh^2 l_k}(C_j+C_i\cosh l_k)(\cosh l_j+\cosh l_i\cosh l_k)],
\end{aligned}
\end{equation}
\begin{equation}\label{Eq: F55}
\begin{aligned}
\frac{\partial \theta_i}{\partial u_j}
=-\frac{1}{A}\frac{1}{\sinh^2 l_k}[&-C_k\sinh^2 l_k +C_j(\cosh l_i+\cosh l_j\cosh l_k)\\
&+C_i(\cosh l_j+\cosh l_i\cosh l_k)],
\end{aligned}
\end{equation}
and
\begin{equation}\label{Eq: F56}
\begin{aligned}
\frac{\partial \theta_i}{\partial u_k}
=-\frac{1}{A}\frac{1}{\sinh^2 l_j}[&-C_j\sinh^2 l_j +C_k(\cosh l_i+\cosh l_j\cosh l_k)\\
&+C_i(\cosh l_k+\cosh l_i\cosh l_j)].
\end{aligned}
\end{equation}

Since $\alpha_i=1$, it follows from (\ref{Eq: F1}) that $e^{f_i}=-\frac{1}{\sinh u_i}$.
As $f_i\rightarrow +\infty$, we deduce $u_i\rightarrow 0^-$.
According to the definition of the admissible space $\mathcal{U}_{ijk}(\eta)$ in (\ref{Eq: admissible space a2}),
neither $u_j$ nor $u_k$ can tend to $-\infty$.
This implies that neither $f_j$ nor $f_k$ can tend to $-\infty$.
Consequently, it suffices to verify that the formula (\ref{Eq: F53}) holds for any constants $a,b,c,d$ under the following four cases:
\begin{description}
\item[(i)] $\lim f_i=+\infty,\ \lim f_j=+\infty,\ \lim f_k=+\infty$,
\item[(ii)] $\lim f_i=+\infty,\ \lim f_j=+\infty,\ \lim f_k=d$,
\item[(iii)] $\lim f_i=+\infty,\ \lim f_j=c,\ \lim f_k=+\infty$,
\item[(iv)] $\lim f_i=+\infty,\ \lim f_j=a,\ \lim f_k=b$.
\end{description}

For simplicity, the coefficients in the cases (i), (ii), (iii), and (iv) are denoted by $a_i, a_j, a_k$ and $b_i, b_j, b_k$, all of which are positive constants.
Note that the formula (\ref{Eq: DCS3}) can be rewritten as
\begin{equation}\label{Eq: F57}
\begin{aligned}
\cosh l_i=&\bigg[\eta_{jk}-\sqrt{(e^{-2f_j}+\alpha_j)
(e^{-2f_k}+\alpha_k)}\bigg]e^{f_j+f_k},\\
\cosh l_j=&\bigg[\eta_{ik}-\sqrt{(e^{-2f_i}+\alpha_i)
(e^{-2f_k}+\alpha_k)}\bigg]e^{f_i+f_k},\\
\cosh l_k=&\bigg[\eta_{ij}-\sqrt{(e^{-2f_i}+\alpha_i)
(e^{-2f_j}+\alpha_j)}\bigg]e^{f_i+f_j}.
\end{aligned}
\end{equation}

For the case (i), if $\lim f_i=+\infty$, $\lim f_j=+\infty$, and $\lim f_k=+\infty$,
then it follows from (\ref{Eq: F57}) that $\lim\cosh l_i:=\lim a_ie^{f_j+f_k}$, $\lim\cosh l_j:=\lim a_je^{f_i+f_k}$, and $\lim\cosh l_k:=\lim a_ke^{f_i+f_j}$.
Moreover, $\lim C_i=\lim  e^{f_i}$,
\begin{eqnarray*}
\lim C_j=
\begin{cases}
1,\ & \text{if}\ \alpha_j=0,  \\
\lim  e^{f_j},\ & \text{if}\ \alpha_j=1,
\end{cases}
\quad \text{and} \quad
\lim C_k=
\begin{cases}
1,\ & \text{if}\ \alpha_k=0,  \\
\lim  e^{f_k},\ & \text{if}\ \alpha_k=1.
\end{cases}
\end{eqnarray*}
Hence, the formula (\ref{Eq: F54}) gives
\begin{equation*}
\begin{aligned}
\lim(-A\frac{\partial \theta_i}{\partial u_i})
=&\, \lim\bigg[\frac{1}{a_j^2e^{2f_i+2f_k}}
(C_k+C_ia_je^{f_i+f_k})(a_ke^{f_i+f_j}+a_ia_je^{f_i+f_j+2f_k})\\
&\ \ \ \ \ \ \ \ \ +\frac{1}{a_k^2e^{2f_i+2f_j}}
(C_j+C_ia_ke^{f_i+f_j})(a_je^{f_i+f_k}+a_ia_ke^{f_i+2f_j+f_k})\bigg]\\
:=&\, \lim b_ie^{f_i+f_j+f_k};
\end{aligned}
\end{equation*}
the formula (\ref{Eq: F55}) gives
\begin{equation*}
\begin{aligned}
\lim(-A\frac{\partial \theta_i}{\partial u_j})
=&\, \lim\bigg[\frac{1}{a_k^2e^{2f_i+2f_j}}
(-C_ka_k^2e^{2f_i+2f_j}
+C_j(a_ie^{f_j+f_k}+a_ja_ke^{2f_i+f_j+f_k})\\
&\ \ \ \ \ \ \ \ \  + C_i(a_je^{f_i+f_k}+a_ia_ke^{f_i+2f_j+f_k}))\bigg]\\
:=&\, b_j e^{f_k};
\end{aligned}
\end{equation*}
and the formula (\ref{Eq: F56}) gives
\begin{equation*}
\begin{aligned}
\lim(-A\frac{\partial \theta_i}{\partial u_k})
=&\, \lim\bigg[\frac{1}{a_j^2e^{2f_i+2f_k}}
(-C_ja_j^2e^{2f_i+2f_k}
+C_k(a_ie^{f_j+f_k}+a_ja_ke^{2f_i+f_j+f_k})\\
& \ \ \ \ \ \ \ \ \ +
C_i(a_ke^{f_i+f_j}+a_ia_je^{f_i+f_j+2f_k}))\bigg]\\
:=&\, \lim b_k e^{f_j}.
\end{aligned}
\end{equation*}
Since
\begin{equation*}
\begin{aligned}
A&=\sinh l_r\sinh l_s\sinh \theta_t\\
&=\sinh l_r\sinh l_s\sqrt{\cosh^2 \theta_t-1}\\
&=\sqrt{\sinh^2 l_r\sinh^2 l_s\cosh^2 \theta_t-\sinh^2 l_r\sinh^2 l_s}\\
&=\sqrt{(\cosh l_t+\cosh l_r\cosh l_s)^2-\sinh^2 l_r\sinh^2 l_s}\\
&=\sqrt{\cosh^2 l_t+2\cosh l_r\cosh l_s\cosh l_t+\cosh(l_r+l_s)\cosh(l_r-l_s)}\\
&>2,
\end{aligned}
\end{equation*}
it follows that $\lim A=\lim \sinh l_r\sinh l_s\sinh \theta_t>0$.
For any constant $C\in \mathbb{R}$, there exists a constant $W>0$, depending on $C$, such that if $f_i\geq W$, then (\ref{Eq: F53}) holds.

For the case $\mathrm{(ii)}$, if $\lim f_i=+\infty,\ \lim f_j=+\infty,\ \lim f_k=d$, then it follows from (\ref{Eq: F57}) that $\lim\cosh l_i:=\lim a_ie^{f_j}$, $\lim\cosh l_j:=\lim a_je^{f_i}$, and $\lim\cosh l_k:=\lim a_ke^{f_i+f_j}$.
Furthermore,
\begin{equation*}
\begin{aligned}
\lim(-A\frac{\partial \theta_i}{\partial u_i})
:=\lim b_ie^{f_i+f_j},\
\lim(-A\frac{\partial \theta_i}{\partial u_j})
:=b_j,\
\lim(-A\frac{\partial \theta_i}{\partial u_k})
:=\lim b_ke^{f_j}.
\end{aligned}
\end{equation*}
Note that $\lim A>0$, for any $C\in \mathbb{R}$, there exists a constant $W>0$, depending on $C$, such that if $f_i\geq W$, then (\ref{Eq: F53}) holds.

For the case (iii), if $\lim f_i=+\infty,\ \lim f_j=c,\ \lim f_k=+\infty$, then it follows from (\ref{Eq: F57}) that $\lim\cosh l_i:=\lim a_ie^{f_k}$, $\lim\cosh l_j:=\lim a_je^{f_i+f_k}$, and $\lim\cosh l_k:=\lim a_ke^{f_i}$.
Furthermore,
\begin{equation*}
\begin{aligned}
\lim(-A\frac{\partial \theta_i}{\partial u_i})
:=\lim b_ie^{f_i+f_k},\
\lim(-A\frac{\partial \theta_i}{\partial u_j})
:=b_je^{f_k},\
\lim(-A\frac{\partial \theta_i}{\partial u_k})
:=\lim b_k.
\end{aligned}
\end{equation*}
Note that $\lim A>0$, for any $C\in \mathbb{R}$, there exists a constant $W>0$, depending on $C$, such that if $f_i\geq W$, then (\ref{Eq: F53}) holds.

For the case (iv), if $\lim f_i=+\infty,\ \lim f_j=a,\ \lim f_k=b$, then it follows from (\ref{Eq: F57}) that $\lim\cosh l_i:=a_i$, $\lim\cosh l_j:=\lim a_je^{f_i}$, and $\lim\cosh l_k:=\lim a_ke^{f_i}$.
Furthermore,
\begin{equation*}
\begin{aligned}
\lim(-A\frac{\partial \theta_i}{\partial u_i})
:=\lim b_ie^{f_i},\
\lim(-A\frac{\partial \theta_i}{\partial u_j})
:=b_j,\
\lim(-A\frac{\partial \theta_i}{\partial u_k})
:=b_k.
\end{aligned}
\end{equation*}
Note that $\lim A>0$, for any $C\in \mathbb{R}$, there exists a constant $W>0$, depending on $C$, such that if $f_i\geq W$, then (\ref{Eq: F53}) holds.
\qed

As a direct corollary of Lemma \ref{lem: CCF a1}, we obtain the following result.

\begin{corollary}\label{Cor: CCF a1}
Under the same assumptions as those in Definition \ref{Def: CCF1},
let $\{ijk\}\in F$ be a hyperbolic right-angled hexagon with $\alpha_i=1$.
For any constants $C_1,C_2,...,C_N\in \mathbb{R}$, there exists a constant $W=W(C_1,C_2,...,C_N)>0$ such that if $f_i\geq W$, then
\begin{equation*}
-\sum_{\{ijk\}\in F}\frac{\partial \theta^{jk}_i}{\partial u_i}
>\sum^N_{j=1,\, j\sim i}C_j\frac{\partial \theta^{jk}_i}{\partial u_j}.
\end{equation*}
\end{corollary}
\proof
By Lemma \ref{lem: CCF a1}, for any constant $C=|C_j|>0$, there exists a constant $W_j=W_j(C_j)>0$ such that if $f_i\geq W_j$, then
\begin{equation*}
-\frac{\partial \theta^{jk}_i}{\partial u_i}
=\bigg|\frac{\partial \theta^{jk}_i}{\partial u_i}\bigg|
>C\bigg(\bigg|\frac{\partial \theta^{jk}_i}{\partial u_j}\bigg|+\bigg|\frac{\partial \theta^{jk}_i}{\partial u_k}\bigg|\bigg)
\geq C_j\frac{\partial \theta^{jk}_i}{\partial u_j},
\end{equation*}
where the first equality follows from $\frac{\partial \theta^{jk}_i}{\partial u_i}<0$ by (\ref{Eq: F54}).
For any constants $C_1,C_2,...,C_N\in \mathbb{R}$, we choose $W=\max_{j\in\{1,2,...,N\}}W_j>0$.
If $f_i\geq W$, then
\begin{equation*}
-\sum_{\{ijk\}\in F}\frac{\partial \theta^{jk}_i}{\partial u_i}
>\sum_{\{ijk\}\in F}C_j\frac{\partial \theta^{jk}_i}{\partial u_j}
=\sum^N_{j=1,\, j\sim i}C_j\frac{\partial \theta^{jk}_i}{\partial u_j}.
\end{equation*}
\qed

\begin{lemma}\label{Lem: boundary a5}
Under the same assumptions as those in Definition \ref{Def: CCF1},
if $\alpha_i=1$, then the solution $u_i(t)$ to the combinatorial Calabi flow (\ref{Eq: CCF1}) is uniformly bounded from above in $\mathbb{R}_{<0}$.
Consequently, it cannot reach the boundary $\partial_{0}\mathcal{U}(\eta)$.
\end{lemma}
\proof
By the arguments in Lemma \ref{Lem: boundary a2},
we assume that there exists $i\in B$ such that $u_i\rightarrow 0^-$.
Thus $f_i\rightarrow+\infty$ and then $\theta^{jk}_{i}\rightarrow 0$ uniformly.
This implies $K_i\rightarrow 0$ uniformly as $f_i\rightarrow+\infty$.
Consequently, there exists constants $W_1>0$ and $W_2>0$, such that if $f_i\geq W_1$, then $K_i-\overline{K}_i<0$ and $|K_i-\overline{K}_i|\geq\frac{1}{2}\overline{K}_i
\geq\frac{1}{2}W_2>0$.
Direct calculations yield
\begin{equation}
\begin{aligned}
\frac{d u_i}{dt}
&=-\Delta(K-\overline{K})_i\\
&=-\frac{\partial K_i}{\partial u_i}(K_i-\overline{K}_i)-\sum_{j\neq i}\frac{\partial K_i}{\partial u_j}(K_j-\overline{K}_j)\\
&=-\sum_{\{ijk\}\in F}\frac{\partial \theta_i^{jk}}{\partial u_i}(K_i-\overline{K}_i)-\sum_{j\sim i}\left(\frac{\partial \theta^{jk}_i}{\partial u_j}+\frac{\partial \theta^{jl}_i}{\partial u_j}\right)(K_j-\overline{K}_j)\\
&=\bigg[-\sum_{\{ijk\}\in F}\frac{\partial \theta_i^{jk}}{\partial u_i}-\sum_{j\sim i}\left(\frac{\partial \theta^{jk}_i}{\partial u_j}+\frac{\partial \theta^{jl}_i}{\partial u_j}\right)\frac{K_j-\overline{K}_j}{K_i-\overline{K}_i}\bigg](K_i-\overline{K}_i).
\end{aligned}
\end{equation}
By Proposition \ref{Prop: property} (ii), the function $\mathcal{C}(u)$ is decreasing along the combinatorial Calabi flow (\ref{Eq: CCF1}).
This implies that $|K_j-\overline{K}_j|$ is bounded for all $j\in B$.
In other words, there exists a constant $W_3>0$, such that $|K_j-\overline{K}_j|\leq W_3$ for any $j\in B$.
Therefore, $\big|\frac{K_j-\overline{K}_j}{K_i-\overline{K}_i}\big|
\leq\frac{2W_3}{W_2}$ for all $j\neq i$.
By Corollary \ref{Cor: CCF a1}, one can choose constants $C_j=C_j(W_2,W_3)$, $j=1,2,...,N$, specially $C_j=\frac{4W_3}{W_2}$ if $\frac{\partial\theta^{jk}_i}{\partial u_j}\geq0$ and $C_j=-\frac{4W_3}{W_2}$ if $\frac{\partial\theta^{jk}_i}{\partial u_j}<0$,
such that there exists $W_4=W_4(C_1,C_2,...,C_N)>0$, if $f_i\geq \max\{W_1,W_4\}$, then
\begin{equation*}
-\sum_{\{ijk\}\in F}\frac{\partial \theta_i^{jk}}{\partial u_i}-\sum_{j\sim i}\left(\frac{\partial \theta^{jk}_i}{\partial u_j}+\frac{\partial \theta^{jl}_i}{\partial u_j}\right)\frac{K_j-\overline{K}_j}{K_i-\overline{K}_i}
>-\sum_{\{ijk\}\in F}\frac{\partial \theta_i^{jk}}{\partial u_i}-\sum^N_{j=1,j\sim i}C_j\frac{\partial \theta^{jk}_i}{\partial u_j}
>0,
\end{equation*}
which implies $\frac{du_i}{dt}<0$.
The remaining part of the proof follows analogously to Lemma \ref{Lem: boundary a2} and is omitted here.
\qed

As an immediate corollary of Lemma \ref{Lem: boundary a4} and Lemma \ref{Lem: boundary a5},
we conclude that the solution $u(t)$ to the combinatorial Calabi flow \eqref{Eq: CCF1} exists for all time,
which corresponds to the existence part of Theorem \ref{Thm: converge} (i).

\begin{corollary}\label{Cor: exist a2}
Under the same assumptions as those in Definition \ref{Def: CCF1}, the solution $u(t)$ to the combinatorial Calabi flow (\ref{Eq: CCF1}) stays in a compact subset of $\mathcal{U}(\eta)$.
Consequently, this solution exists for all time.
\end{corollary}

The following theorem shows the convergence of the solution $u(t)$ to the combinatorial Calabi flow (\ref{Eq: CCF1}),
which corresponds to the convergence part of Theorem \ref{Thm: converge} (i).

\begin{theorem}\label{Thm: converge a2}
Under the same assumptions as those in Definition \ref{Def: CCF1},
the solution $u(t)$ to the combinatorial Calabi flow (\ref{Eq: CCF1}) converges exponentially fast.
\end{theorem}
\proof
By combining Corollary \ref{Cor: exist a2}, Theorem \ref{Thm: matrix negative 1}, and the continuity of the eigenvalues of $\Delta$,
there exists a constant $\lambda_0>0$ such that all eigenvalues $\lambda_\Delta$ of $\Delta$ satisfy $\lambda_\Delta<-\sqrt{\lambda_0/2}$ along the combinatorial Calabi flow (\ref{Eq: CCF1}).
Therefore, along the combinatorial Calabi flow (\ref{Eq: CCF1}),
\begin{equation*}
\frac{d\mathcal{C}(u(t))}{dt}
=-(K-\overline{K})^T\Delta^2(K-\overline{K})\leq -\lambda_0\mathcal{C}(u(t)),
\end{equation*}
which implies
\begin{equation*}
\mathcal{C}(u(t))
=\frac{1}{2}\Vert K(t)-\overline{K}\Vert^2
\leq e^{-\lambda_0t}\Vert K(0)-\overline{K}\Vert^2.
\end{equation*}
By combining Theorem \ref{Thm: rigidity and image} (i) and Corollary \ref{Cor: exist a2}, it follows that
\begin{equation*}
\Vert u(t)-\overline{u}\Vert^2
\leq C_1\Vert K(t)-\overline{K}\Vert^2
\leq 2C_1e^{-\lambda_0t}\Vert K(0)-\overline{K}\Vert^2
\leq C_2e^{-\lambda_0t}
\end{equation*}
for some positive constants $C_1,C_2$.
This completes the proof.
\qed

\begin{remark}
The exponential convergence of the solution $u(t)$ to the combinatorial Calabi flow (\ref{Eq: CCF1}) can also be established by Lyapunov Stability Theorem (\cite[Chapter 5]{Pontryagin}), following an argument analogous to that used in the proof of Theorem \ref{Thm: converge a1}.
\end{remark}

\section{Combinatorial curvature flows for the discrete conformal structure (\ref{Eq: DCS1})}\label{Sec: DSC1}

Suppose $(\Sigma,\mathcal{T},\eta)$ is a weighted triangulated surface with boundary, where $\eta\in \mathbb{R}_{>0}^E$.
According to Lemma 2.21 and Remark 2.22 in \cite{X-Z DCS2}, for any $i\in B$,
\begin{equation}\label{Eq: F2}
u_i=-e^{-f_i}.
\end{equation}
Furthermore, we have the following results.

\begin{theorem}[\cite{X-Z DCS2}]
\label{Thm: ASC3}
Suppose $(\Sigma,\mathcal{T},\eta)$ is a weighted triangulated surface with boundary, where $\eta\in \mathbb{R}_{>0}^E$.
The admissible space
\begin{equation}\label{Eq: F3}
\mathcal{U}_{ijk}(\eta)
=\{(u_i,u_j,u_k)\in \mathbb{R}_{<0}^3 \mid u_r+u_s>-\sqrt{2\eta_{rs}}, \ \forall\{r,s\}\subseteq\{i,j,k\} \},
\end{equation}
is a convex polytope.
As a result, the admissible space
$\mathcal{U}(\eta)=\bigcap_{\{ijk\}\in F}\mathcal{U}_{ijk}(\eta)$
is a convex polytope.
\end{theorem}

\begin{theorem}[\cite{X-Z DCS2}, Theorem 5.2]
\label{Thm: matrix negative 2}
Suppose $(\Sigma,\mathcal{T},\eta)$ is a weighted triangulated surface with boundary, where $\eta\in \mathbb{R}_{>0}^E$.
For a right-angled hyperbolic hexagon $\{ijk\}\in F$ on $\mathcal{U}_{ijk}(\eta)$,
the Jacobian $\Lambda_{ijk}=\frac{\partial(\theta^{jk}_i,\theta^{ik}_j,\theta^{ij}_k)}
{\partial(u_i,u_j,u_k)}$ is symmetric and negative definite.
Consequently, the Jacobian $\Lambda=\frac{\partial (K_i,..., K_N)}{\partial(u_i,...,u_N)}$ is symmetric and negative definite on $\mathcal{U}(\eta)$.
\end{theorem}

The combinatorial curvature flows for the discrete conformal structure (\ref{Eq: DCS1}) is redefined as follows.

\begin{definition}\label{Def: CCF2}
Suppose $(\Sigma,\mathcal{T},\eta)$ is a weighted triangulated surface with boundary, where $\eta\in \mathbb{R}_{>0}^E$.
Let $\overline{K}\in \mathbb{R}_{>0}^N$ be a given function defined on $B=\{1,2,...,N\}$.
The combinatorial Ricci flow for the discrete conformal structure (\ref{Eq: DCS1}) is defined as
\begin{eqnarray}\label{Eq: CRF2}
\begin{cases}
\frac{du_i}{dt}=K_i-\overline{K}_i,\\
u_i(0)=u_0.
\end{cases}
\end{eqnarray}
The combinatorial Calabi flow for the discrete conformal structure (\ref{Eq: DCS1}) is defined as
\begin{eqnarray}\label{Eq: CCF2}
\begin{cases}
\frac{du_i}{dt}=-\Delta(K-\overline{K})_i,\\
u_i(0)=u_0,
\end{cases}
\end{eqnarray}
where $\Delta=(\frac{\partial K_i}{\partial u_j})_{N\times N}$ is the discrete Laplace operator in Theorem \ref{Thm: matrix negative 2}.
\end{definition}

By (\ref{Eq: F3}), the boundary of the admissible space $\mathcal{U}(\eta)$ in $[-\infty,0]^{N}$ consists of the following three parts:
\begin{description}
\item[(i)]
$\partial_{\infty}\mathcal{U}(\eta)=\{u\in[-\infty,0]^N \mid \text{there exists}\, i\in B\, \text{such that}\,  u_i=-\infty \}$,
\item[(ii)]
$\partial_0\mathcal{U}(\eta)=\{u\in[-\infty,0]^N \mid \text{there exists}\, i\in B\, \text{such that}\,  u_i=0 \}$,
\item[(iii)]
$\partial_{l}\mathcal{U}(\eta)
=\bigcup_{\{ij\}\in E}\partial_{ij}\mathcal{U}(\eta)
=\bigcup_{\{ij\}\in E}\{u\in [-\infty,0]^N \mid u_i+u_j=-\sqrt{2\eta_{ij}}\}$.
\end{description}

\subsection{Combinatorial Ricci flow for the discrete conformal structure (\ref{Eq: DCS1})}

The following lemma shows that the solution $u(t)$ to the combinatorial Ricci flow (\ref{Eq: CRF2}) cannot reach the boundaries $\partial_{\infty}\mathcal{U}(\eta)$ and $\partial_{l}\mathcal{U}(\eta)$.
Since the proof is analogous to those of Lemma \ref{Lem: boundary a1} and Lemma \ref{Lem: boundary a3}, it is omitted here.

\begin{lemma}\label{Lem: boundary b1}
Suppose $(\Sigma,\mathcal{T},\eta)$ is a weighted triangulated surface with boundary, where $\eta\in \mathbb{R}_{>0}^E$.
The solution $u(t)$ to the combinatorial Ricci flow (\ref{Eq: CRF2}) stays in a bounded subset of $\mathbb{R}_{<0}^N$ and cannot reach the boundary $\partial_{l}\mathcal{U}(\eta)$.
\end{lemma}

To prove that the solution $u(t)$ to the combinatorial Ricci flow (\ref{Eq: CRF2}) cannot reach the boundary $\partial_{0}\mathcal{U}(\eta)$, we need the following lemma.

\begin{lemma}(\cite{GL2}, Lemma 4.6)\label{Lem: GL-limit}
Suppose $\{ijk\}\in F$ is a right-angled hyperbolic  hexagon adjacent to the boundary components $i,j,k\in B$.
Then
\begin{equation*}
\lim_{f_k\rightarrow+\infty}\theta^{ij}_k(f_i,f_j,f_k)=0
\end{equation*}
and the converge is uniform.
\end{lemma}

\begin{lemma}\label{Lem: boundary b2}
Suppose $(\Sigma,\mathcal{T},\eta)$ is a weighted triangulated surface with boundary, where $\eta\in \mathbb{R}_{>0}^E$.
The solution $u_i(t)$ to the combinatorial Ricci flow (\ref{Eq: CRF2}) is uniformly bounded from above in $\mathbb{R}_{<0}$.
Consequently, it cannot reach the boundary $\partial_{0}\mathcal{U}(\eta)$.
\end{lemma}
\proof
If $\lim_{t\rightarrow T}u_i(t)=0^-$ for $T\in(0,+\infty]$, then $\lim_{t\rightarrow T}f_i(t)=+\infty$ by (\ref{Eq: F2}).
By Lemma \ref{Lem: GL-limit}, we have
$\theta^{jk}_{i}\rightarrow 0^-$ uniformly as $f_i\rightarrow+\infty$.
The remainder of the proof is analogous to that of Lemma \ref{Lem: boundary a2}, and thus is omitted here.
\qed

A direct corollary of Lemma \ref{Lem: boundary b1} and Lemma \ref{Lem: boundary b2} shows the longtime existence of the solution $u(t)$ to the combinatorial Ricci flow (\ref{Eq: CRF2}),
which corresponds to the existence part of Theorem \ref{Thm: converge} (ii).

\begin{corollary}\label{Cor: exist b1}
Suppose $(\Sigma,\mathcal{T},\eta)$ is a weighted triangulated surface with boundary, where $\eta\in \mathbb{R}_{>0}^E$.
The solution $u(t)$ to the combinatorial Ricci flow (\ref{Eq: CRF2}) stays in a compact subset of $\mathcal{U}(\eta)$.
As a result, this solution exists for all time.
\end{corollary}

The following theorem shows the convergence of the solution $u(t)$ to the combinatorial Ricci flow (\ref{Eq: CRF2}),
which corresponds to the convergence part of Theorem \ref{Thm: converge} (ii).
Since the proof is analogous to that of Theorem \ref{Thm: converge a1}, it is omitted here.

\begin{theorem}\label{Thm: converge b1}
Suppose $(\Sigma,\mathcal{T},\eta)$ is a weighted triangulated surface with boundary, where $\eta\in \mathbb{R}_{>0}^E$.
The solution to the combinatorial Ricci flow (\ref{Eq: CRF2}) converges exponentially fast.
\end{theorem}

\subsection{Combinatorial Calabi flow for the discrete conformal structure (\ref{Eq: DCS1})}

The following lemma shows that the solution $u(t)$ to the combinatorial Calabi flow (\ref{Eq: CCF2}) cannot reach the boundaries $\partial_{\infty}\mathcal{U}(\eta)$ and $\partial_{l}\mathcal{U}(\eta)$.
The proof is identical and thus omitted.

\begin{lemma}\label{Lem: boundary b4}
Suppose $(\Sigma,\mathcal{T},\eta)$ is a weighted triangulated surface with boundary, where $\eta\in \mathbb{R}_{>0}^E$.
The solution $u(t)$ to the combinatorial Calabi flow (\ref{Eq: CCF2}) stays in a bounded subset of $\mathbb{R}_{<0}^N$ and cannot reach the boundary $\partial_{l}\mathcal{U}(\eta)$.
\end{lemma}

To prove the solution $u(t)$ to the combinatorial Calabi flow (\ref{Eq: CCF2}) cannot reach the boundary $\partial_{0}\mathcal{U}(\eta)$, we need the following lemma.

\begin{lemma}\label{lem: CCF b1}
Suppose $(\Sigma,\mathcal{T},\eta)$ is a weighted triangulated surface with boundary, where $\eta\in \mathbb{R}_{>0}^E$.
Then for any constant $C\in \mathbb{R}$, there exists a constant $W=W(C)>0$ such that if $f_i\geq W$, then
\begin{equation}\label{Eq: F58}
\bigg|\frac{\partial \theta^{jk}_i}{\partial u_i}\bigg|
>C\bigg(\bigg|\frac{\partial \theta^{jk}_i}{\partial u_j}\bigg|+\bigg|\frac{\partial \theta^{jk}_i}{\partial u_k}\bigg|\bigg).
\end{equation}
\end{lemma}
\proof
For simplicity, set $l_r=l_{st}$, $\theta_r=\theta^{st}_r$, and $A=\sinh l_r\sinh l_s\sinh \theta_t$, where $\{r,s,t\}=\{i,j,k\}$.
Replacing $C_r$ in (\ref{Eq: F4}) with $e^{f_r}$ yields the matrix $\frac{\partial(\theta^{jk}_i,\theta^{ik}_j,\theta^{ij}_k)}
{\partial(u_i,u_j,u_k)}$ in Theorem \ref{Thm: matrix negative 2},
a detailed derivation can be found in Theorem 5.2 of \cite{X-Z DCS2}.
Moreover, it follows from (\ref{Eq: d1}) that
\begin{equation*}
\begin{aligned}
\coth d_{ij}
&=\frac{1}{\sinh l_k}[\sinh(f_j-f_i)+\eta_{ij}e^{f_i+f_j}]\\
&=\frac{1}{\sinh l_k}[\sinh(f_j-f_i)+\cosh(f_j-f_i)+\cosh l_k]\\
&=\frac{e^{f_j}+\cosh l_ke^{f_i}}{\sinh l_ke^{f_i}},
\end{aligned}
\end{equation*}
which aligns with (\ref{Eq: F5}) under the substitution $C_r$ with $e^{f_r}$.
Therefore, replacing $C$ in (\ref{Eq: F54}), (\ref{Eq: F55}) and (\ref{Eq: F56}) with $e^f$ yields
\begin{equation*}
\begin{aligned}
\frac{\partial \theta_i}{\partial u_i}
=-\frac{1}{A}[&\frac{1}{\sinh^2 l_j}
(e^{f_k}+e^{f_i}\cosh l_j)(\cosh l_k+\cosh l_i\cosh l_j)\\
&+\frac{1}{\sinh^2 l_k}
(e^{f_j}+e^{f_i}\cosh l_k)(\cosh l_j+\cosh l_i\cosh l_k)],\\
\frac{\partial \theta_i}{\partial u_j}
=-\frac{1}{A}&\frac{1}{\sinh^2 l_k}[-e^{f_k}\sinh^2 l_k +e^{f_j}(\cosh l_i+\cosh l_j\cosh l_k)\\
&\ \ \ \ \ \ \ \ \ \ \ \ \ \ +e^{f_i}(\cosh l_j+\cosh l_i\cosh l_k)],\\
\frac{\partial \theta_i}{\partial u_k}
=-\frac{1}{A}&\frac{1}{\sinh^2 l_j}[-e^{f_j}\sinh^2 l_j +e^{f_k}(\cosh l_i+\cosh l_j\cosh l_k)\\
&\ \ \ \ \ \ \ \ \ \ \ \ \ \ +e^{f_i}(\cosh l_k+\cosh l_i\cosh l_j)].
\end{aligned}
\end{equation*}

By Lemma \ref{Lem: boundary b4}, the solution $u(t)$ to the combinatorial Calabi flow (\ref{Eq: CCF2}) stays in a bounded subset of $\mathbb{R}_{<0}^N$.
This implies that $f$ cannot tend to $-\infty$ by (\ref{Eq: F2}).
Consequently, it suffices to verify that the formula (\ref{Eq: F58}) holds for any constants $a,b,c$ under the following three cases:
\begin{description}
\item[(i)] $\lim f_i=+\infty,\ \lim f_j=+\infty,\ \lim f_k=+\infty$,
\item[(ii)] $\lim f_i=+\infty,\ \lim f_j=+\infty,\ \lim f_k=c$,
\item[(iii)] $\lim f_i=+\infty,\ \lim f_j=a,\ \lim f_k=b$.
\end{description}

For simplicity, the coefficients in the cases (i), (ii), and (iii) are denoted by $a_i, a_j, a_k$ and $b_i, b_j, b_k$, all of which are positive constants.
Note that the formula (\ref{Eq: DCS1}) can be reformulated as
\begin{gather*}
\cosh l_k=-\cosh (f_j-f_i)+\eta_{ij}e^{f_i+f_j}
=(\eta_{ij}-\frac{1}{2}e^{-2f_i}-\frac{1}{2}e^{-2f_j})e^{f_i+f_j},\\
\cosh l_j=-\cosh (f_k-f_i)+\eta_{ik}e^{f_i+f_k}
=(\eta_{ik}-\frac{1}{2}e^{-2f_i}-\frac{1}{2}e^{-2f_k})e^{f_i+f_k},\\
\cosh l_i=-\cosh (f_j-f_k)+\eta_{jk}e^{f_j+f_k}
=(\eta_{jk}-\frac{1}{2}e^{-2f_j}-\frac{1}{2}e^{-2f_k})e^{f_j+f_k}.
\end{gather*}

For the case (i), if $\lim f_i=+\infty,\ \lim f_j=+\infty,\ \lim f_k=+\infty$, then
$\lim \cosh l_i=\lim a_ie^{f_j+f_k}$,\
$\lim \cosh l_j=\lim a_je^{f_i+f_k}$ and
$\lim \cosh l_k=\lim a_ke^{f_i+f_j}$.
Hence,
\begin{align*}
\lim(-A\frac{\partial \theta_i}{\partial u_i})
&=\lim\bigg[\frac{1}{a^2_je^{2f_i+2f_k}}
(e^{f_k}+a_je^{2f_i+f_k})(a_ke^{f_i+f_j}+a_ia_je^{f_i+f_j+2f_k})
\\
&\ \ \ \ \ \ \ \ \ \ \ \ \  +\frac{1}{a^2_ke^{2f_i+2f_j}}
(e^{f_j}+a_ke^{2f_i+f_j})(a_je^{f_i+f_k}+a_ia_ke^{f_i+2f_j+f_k})
\bigg]\\
&:=\lim b_ie^{f_i+f_j+f_k};\\
\lim(-A\frac{\partial \theta_i}{\partial u_j})
=&\lim\bigg[\frac{1}{a^2_ke^{2f_i+2f_j}}
(-a^2_ke^{2f_i+2f_j+f_k}+e^{f_j}
(a_ie^{f_j+f_k}+a_ja_ke^{2f_i+f_j+f_k})\\
&\ \ \ \ \ \ \ \ \ \ \ \ \  +e^{f_i}(a_je^{f_i+f_k}+a_ia_ke^{f_i+2f_j+f_k}))\bigg]\\
&:=\lim b_je^{f_k};\\
\lim(-A\frac{\partial \theta_i}{\partial u_k})
&=\lim\bigg[\frac{1}{a^2_je^{2f_i+2f_k}}
(-a^2_je^{2f_i+f_j+2f_k}+e^{f_k}
(a_ie^{f_j+f_k}+a_ja_ke^{2f_i+f_j+f_k})\\
&\ \ \ \ \ \ \ \ \ \ \ \ \  +e^{f_i}(a_ke^{f_i+f_j}+a_ia_je^{f_i+f_j+2f_k}))\bigg]\\
&:=\lim b_ke^{f_j}.
\end{align*}
Given $\lim A>0$, for any $C\in \mathbb{R}$, there exists a constant $W>0$, depending on $C$, such that if $f_i\geq W$, then (\ref{Eq: F58}) holds.
The proofs for cases (ii) and (iii) follow analogous reasoning and are omitted here.
\qed

As a direct corollary of Lemma \ref{lem: CCF b1}, we obtain the following result.

\begin{corollary}\label{Cor: CCF b1}
Suppose $(\Sigma,\mathcal{T},\eta)$ is a weighted triangulated surface with boundary, where $\eta\in \mathbb{R}_{>0}^E$.
For any constants $C_1,C_2,...,C_N\in \mathbb{R}$, there exists a constant $W=W(C_1,C_2,...,C_N)>0$ such that if $f_i\geq W$, then
\begin{equation*}
-\sum_{\{ijk\}\in F}\frac{\partial \theta^{jk}_i}{\partial u_i}
>\sum^N_{j=1,j\sim i}C_j\frac{\partial \theta^{jk}_i}{\partial u_j}.
\end{equation*}
\end{corollary}

The proof of the following lemma is analogous to that of Lemma \ref{Lem: boundary a5} and thus is omitted.

\begin{lemma}\label{Lem: boundary b5}
Suppose $(\Sigma,\mathcal{T},\eta)$ is a weighted triangulated surface with boundary, where $\eta\in \mathbb{R}_{>0}^E$.
The solution $u_i(t)$ to the combinatorial Calabi flow (\ref{Eq: CCF2}) is uniformly bounded from above in $\mathbb{R}_{<0}$.
Consequently, it cannot reach the boundary $\partial_{0}\mathcal{U}(\eta)$.
\end{lemma}

A direct corollary of Lemma \ref{Lem: boundary b4} and Lemma \ref{Lem: boundary b5} implies the longtime existence to the solution $u(t)$ of the combinatorial Calabi flow (\ref{Eq: CCF2}),
which corresponds to the existence part of Theorem \ref{Thm: converge} (ii).

\begin{corollary}\label{Cor: exist b2}
Suppose $(\Sigma,\mathcal{T},\eta)$ is a weighted triangulated surface with boundary, where $\eta\in \mathbb{R}_{>0}^E$.
The solution $u(t)$ to the combinatorial Calabi flow (\ref{Eq: CCF2}) stays in a compact subset of the admissible space $\mathcal{U}(\eta)$.
As a result, this solution exists for all time.
\end{corollary}

The following theorem shows the convergence of the solution $u(t)$ to the combinatorial Calabi flow (\ref{Eq: CCF2}),
which corresponds to the convergence part of Theorem \ref{Thm: converge} (ii).
The proof is analogous to that of Theorem \ref{Thm: converge a2} and thus is omitted.

\begin{theorem}\label{Thm: converge b2}
Suppose $(\Sigma,\mathcal{T},\eta)$ is a weighted triangulated surface with boundary, where $\eta\in \mathbb{R}_{>0}^E$.
The solution to the combinatorial Calabi flow (\ref{Eq: CCF2}) converges exponentially fast.
\end{theorem}

\section{Combinatorial curvature flows for the remaining discrete conformal structures}\label{Sec: remaining}

In this section, we consider the combinatorial curvature flows for the remaining discrete conformal structures, namely,
the discrete conformal structure (\ref{Eq: new 1}),
along with the mixed discrete conformal structures I, II, and III.

The following theorem is a restatement of Theorem \ref{Thm: converge} (iii).

\begin{theorem}\label{Thm: converge c1}
For the discrete conformal structure (\ref{Eq: new 1}), and the mixed discrete conformal structures I and III,
let $\overline{K}\in (0,+\infty)^N$ be a given function defined on $B=\{1,2,...,N\}$.
Under the assumptions of Theorem \ref{Thm: rigidity and image} (ii), (iv) and (vi),
there exists a constant $\delta>0$ such that if $\Vert K(u(0))-\overline{K}\Vert<\delta$,
then the solutions to the combinatorial Ricci flow (\ref{Eq: CRF}) and the combinatorial Calabi flow (\ref{Eq: CCF})  exist for all time and converge exponentially fast, respectively.
\end{theorem}
\proof
For the combinatorial Ricci flow (\ref{Eq: CRF}), set $\Gamma(u)=K-\overline{K}$.
Then $D\Gamma|_{u=\overline{u}}=\Delta$ is negative definite,
which implies $\overline{u}$ is a local attractor of the combinatorial Ricci flow (\ref{Eq: CRF}).
Then the conclusion follows from Lyapunov Stability Theorem (\cite[Chapter 5]{Pontryagin}).
Similarly, for the combinatorial Calabi flow (\ref{Eq: CCF}), set $\Gamma(u)=-\Delta(K-\overline{K})$.
Then $D\Gamma|_{u=\overline{u}}=-\Delta^2$ is negative definite,
which implies $\overline{u}$ is a local attractor of the combinatorial Calabi flow (\ref{Eq: CCF}).
Then the conclusion follows from Lyapunov Stability Theorem (\cite[Chapter 5]{Pontryagin}).
\qed

The following theorem is a restatement of Theorem \ref{Thm: converge} (iv).

\begin{theorem}\label{Thm: converge c2}
Suppose $(\Sigma,\mathcal{T},\eta)$ is a weighted triangulated surface with boundary, where $\eta\in [1,+\infty)^E$.
For the mixed discrete conformal structure II, let $\overline{K}\in (0,+\infty)^N$ be a given function defined on $B=\{1,2,...,N\}$.
\begin{description}
  \item[(1)] If the solution $u(t)$ to the combinatorial Ricci flow (\ref{Eq: CRF}) or the combinatorial Calabi flow (\ref{Eq: CCF}) converges to $\overline{u}\in \mathcal{U}(\eta)$, then $K(\overline{u})=\overline{K}$.
  \item[(2)] Suppose that there exists a discrete conformal factor $\overline{u}$ on $(\Sigma,\mathcal{T})$ with $K(\overline{u})=\overline{K}$, then there exists a constant $\delta>0$ such that if $\Vert K(u(0))-\overline{K}\Vert<\delta$, the solutions to the combinatorial Ricci flow (\ref{Eq: CRF}) and the combinatorial Calabi flow (\ref{Eq: CCF})  exist for all time and converge exponentially fast to $\overline{u}$.
\end{description}
\end{theorem}
\proof
\noindent\textbf{(1):} Let $u(t)$ be a solution of the combinatorial Ricci flow (\ref{Eq: CRF}).
If $\overline{u}:=u(+\infty)=\lim_{t\rightarrow +\infty} u(t)$ exists in $\mathcal{U}(\eta)$, then $K(\overline{u})=\lim_{t\rightarrow +\infty}K(u(t))$ exists by the $C^1$-smoothness of $K$.
Furthermore, there exists a sequence $\xi_n\in(n,n+1)$ such that as $n\rightarrow +\infty$,
\begin{equation*}
u_i(n+1)-u_i(n)=u'_i(\xi_n)
=K_i(u(\xi_n))-\overline{K}_i\rightarrow 0,
\end{equation*}
which implies $K(\overline{u})=\overline{K}$.
Similarly, if the solution $u(t)$ of the combinatorial Calabi flow (\ref{Eq: CCF}) converges, then
$K(\overline{u})=\lim_{t\rightarrow +\infty}K(u(t))$ exists by the $C^1$-smoothness of $K$.
Furthermore, there exists a sequence $\xi_n\in(n,n+1)$ such that as $n\rightarrow +\infty$,
\begin{equation*}
u_i(n+1)-u_i(n)=u'_i(\xi_n)
=-\Delta(K(u(\xi_n))-\overline{K})_i\rightarrow 0,
\end{equation*}
which implies $K(\overline{u})=\overline{K}$ by the strictly negative definiteness of $\Delta$.

\noindent\textbf{(2):}
The proof is analogous to that of Theorem \ref{Thm: converge c1} and thus is omitted here.
\qed

\begin{remark}
Theorem \ref{Thm: converge c1} and Theorem \ref{Thm: converge c2} establish the longtime existence and convergence of solutions to the combinatorial Ricci flow and combinatorial Calabi flow for initial data with small energy.
However, for general initial data, we have not yet shown that these solutions to the combinatorial curvature flows cannot reach some boundary of the admissible space,
which prevents us from proving their longtime existence and global convergence on surfaces with boundary.
Specifically, for the discrete conformal structure (\ref{Eq: new 1}), the discrete conformal factors lie in the domain $u\in (-\frac{\pi}{2},0)^N$.
Analogous to the proof of Lemma \ref{Lem: boundary a2}, one can show that solution $u_i(t)$ to the combinatorial Ricci flow (\ref{Eq: CRF}) is uniformly bounded from above in $\mathbb{R}_{<0}$.
A critical unresolved issue, however, is proving that the solution $u_i(t)$ remain uniformly bounded from below away from $-\frac{\pi}{2}$.
For the mixed discrete conformal structures III, where discrete conformal factors lie in $u\in\mathbb{R}_{>0}^{N_1}\times\mathbb{R}_{<0}^{N-N_1}$,
we similarly lack a proof for the uniform lower boundedness of $u_i(t)$ in $\mathbb{R}_{>0}$ for the combinatorial Ricci flow (\ref{Eq: CRF}).
As for mixed discrete conformal structures II, the underlying reason for the absence of such results is analogous to that of the discrete conformal structure (\ref{Eq: new 1}).
For the mixed discrete conformal structures I, in the context of specific cases, one can establish the longtime existence and global convergence of solutions to the combinatorial curvature flows with general initial values.
Given the extraordinary complexity inherent in the analysis, a detailed elaboration is omitted herein.
\end{remark}

\end{document}